\documentclass[titlepage,11pt]{article}
\usepackage{latexsym}
\usepackage{amsfonts}
\usepackage{amsmath}
\usepackage{tikz, fp}
\usepackage{textcomp}
\usetikzlibrary{graphs}
\usepackage{float}
\usetikzlibrary{decorations.markings}
\usetikzlibrary{decorations,decorations.pathmorphing}
\newcommand\blackslug{\hbox{\hskip 1pt \vrule width 4pt height 8pt depth 1.5pt
        \hskip 1pt}}
\newcommand\bbox{\hfill \quad \blackslug \bigbreak}

\def\ll{,\ldots,}
\renewcommand{\epsilon}{\varepsilon}
\newcommand{\vare}{\varepsilon}
\DeclareMathOperator{\polylog}{polylog}

\title{Pure pairs. I. Trees and linear anticomplete pairs}
\author{Maria Chudnovsky\thanks{Supported by NSF grant DMS-1550991.
This material is based upon work supported in part by the U. S. Army
Research Laboratory and the U. S. Army Research Office under grant
number
W911NF-16-1-0404.}\\
Princeton University, Princeton, NJ 08544, USA
\\
\\
Alex Scott\thanks{Supported by a Leverhulme Trust Research
Fellowship.}\\
Mathematical Institute, University of Oxford, Oxford OX2 6GG, UK
\\
\\
Paul Seymour\thanks{Supported by ONR grant N00014-14-1-0084, and NSF
grant DMS-1265563, and AFOSR grant A9550-19-1-0187.}\\
Princeton University, Princeton, NJ 08544, USA
\\
\\
Sophie Spirkl\thanks{This work was mostly performed while Spirkl was at Princeton University.}\\
University of Waterloo, Waterloo, Ontario N2L3G1, Canada.}

\date{March 4, 2018; revised \today}

\newtheorem{thm}{}[section]

\newcommand{\Proof}{\noindent{\bf Proof.}\ \ }

\begin{document}
\maketitle
\begin{abstract}
%We prove a conjecture of Liebenau, Pilipczuk, and the last two authors~\cite{cats}, 
%that for every forest $H$ there exists 
%$\epsilon>0$, such that if $G$ is a graph with $|G|>1$ that does not contain $H$ as an induced subgraph, then 
%either 
%\begin{itemize}
%\item some vertex has degree at least $\epsilon |G|$; or 
%\item there are two disjoint sets $A,B\subseteq V(G)$ with $|A|,|B|\ge \epsilon |G|$, such that there are no edges between $A,B$.
%\end{itemize}
%(It is known that no graphs $H$ except forests have this property.)
%Consequently we prove that for every forest $H$, there exists $c>0$ such that for every graph $G$
%containing neither $H$ nor its complement as an induced subgraph, there is a clique or stable set of cardinality at
%least $|G|^c$. 
The Erd\H os-Hajnal conjecture asserts that for every graph $H$ there is a constant $c>0$ such that every graph $G$ that does not contain $H$ as an induced subgraph has a clique or stable set of cardinality at least $|G|^c$.  
In this paper, we prove a conjecture
of  Liebenau and Pilipczuk~\cite{cats}, that
for every forest $H$ there exists 
$c>0$, such that every graph $G$ with $|G|>1$ contains either an induced copy of $H$, or a vertex of degree at least $c|G|$, or 
two disjoint sets of at least $c|G|$ vertices with no edges between them.
It follows that for every forest $H$ there exists $c>0$ such that, if $G$ contains neither $H$ nor its complement as an induced 
subgraph, then there is a clique or stable set of cardinality at least $|G|^c$.   
\end{abstract}

\section{Introduction}

All graphs in this paper are finite and have no loops
or parallel edges. The number of vertices of $G$ is denoted by $|G|$.
If $G,H$ are graphs, we say $G$ {\em contains} $H$ if some induced subgraph of $G$ is isomorphic to $H$, and $G$ is
{\em $H$-free} otherwise. We denote by $\alpha(G), \omega(G)$ the cardinalities of the largest stable sets and largest cliques in $G$
respectively. Two disjoint sets $A,B$ are {\em complete} if every vertex in $A$ is adjacent to every vertex in $B$, and
{\em anticomplete} if no vertex in $A$ has a neighbour in $B$; and we say $A$ {\em covers} $B$ if
every vertex in $B$ has a neighbour in $A$. A pair $(A,B)$ of subsets of $V(G)$ is {\em pure} if $A\cap B=\emptyset$ and $A$ 
is either complete or anticomplete to $B$. This is the first of a series of papers concerning pure pairs.

The Erd\H{o}s-Hajnal conjecture~\cite{EH0, EH} asserts that:
\begin{thm}\label{EHconj}
{\bf Conjecture: }For every graph $H$, there exists $c>0$ such that every $H$-free graph $G$ satisfies
$$\alpha(G)\omega(G)\ge |G|^c.$$
\end{thm}
One way to try to prove this for appropriate graphs $H$ might be to prove a stronger property.  To discuss this, it is helpful to use the language of graph ideals. 

An {\em ideal} of graphs is a class of graphs closed under isomorphism and under taking induced subgraphs; and an 
ideal is {\em proper} if it is not the class of all graphs.  Let us say that an ideal $\mathcal G$  has the {\em Erd\H{o}s-Hajnal property} if
there is some $\epsilon>0$ such that every graph $G\in\mathcal G$ has a clique or stable set of size at least $|G|^\epsilon$.  (Thus the 
Erd\H os-Hajnal conjecture says that for every graph $H$, the ideal of $H$-free graphs has the Erd\H{o}s-Hajnal property; equivalently, every 
proper ideal has the Erd\H os-Hajnal property.)  
An ideal $\mathcal G$  has the {\em strong Erd\H os-Hajnal property} if
there is some $\epsilon>0$ such that every graph $G\in\mathcal G$ with at least two vertices 
contains disjoint sets $A$, $B$ that have size at least $\epsilon |G|$ and are either complete or anticomplete  (that is, $(A,B)$ is a pure pair).
It is easy to prove that if an ideal has the strong Erd\H{o}s-Hajnal property then it has the Erd\H{o}s-Hajnal property
(see~\cite{APPRS,fp}).  

Unfortunately, not every proper ideal has the strong Erd\H{o}s-Hajnal property: by considering sparse random graphs, it is easy to show 
that if the ideal of $H$-free graphs has the Erd\H os-Hajnal property then $H$ must be a forest; and the same argument in the complement 
shows that the complement of $H$ must also be a forest.  Thus $H$ must have at most four vertices, and 
the conjecture is already resolved in these cases.  

However, what if we exclude more than one graph?  In particular, if we consider ideals defined by a finite number of excluded induced 
subgraphs, then when do we get the strong Erd\H os-Hajnal property?  The random graph argument again shows that one of the excluded 
graphs must be a forest and one of them must be the complement of a forest.  Could this be enough?   For example, Bousquet, Lagoutte and 
Thomass\'e~\cite{lagoutte} proved that
for every path $H$, every graph which is both $H$-free and $\overline{H}$-free (with at least two vertices) 
does have a pure pair of sets both of linear size. 
This was extended by Choromanski, Falik, Liebenau, Patel, and Pilipczuk~\cite{hooks}, who proved the same when $H$
is a path with a leaf added adjacent to the third vertex; and then extended 
further by  Liebenau, Pilipczuk, Seymour and Spirkl~\cite{cats}, who proved that ``subdivided caterpillars'' have the property.
Liebenau and Pilipczuk (published in~\cite{cats}) conjectured that in fact all forests have the property:
\begin{thm}\label{forestsymm}
For every forest $H$, there exists $\epsilon>0$ such that for every graph $G$ with $|G|>1$ that is both $H$-free and 
$\overline{H}$-free, there is a pure pair $(A,B)$ with $|A|,|B|\ge \epsilon |G|$.
\end{thm}
This is true, and a consequence of our main result.  In light of the random graph argument, this completes the classification of ideals that are defined by a finite number of excluded induced subgraphs and have the strong Erd\H os-Hajnal property.

Since the strong Erd\H os-Hajnal property implies the Erd\H os-Hajnal property,
it follows from \ref{forestsymm} that
\begin{thm}\label{treeEH}
For every forest $H$, there exists $c>0$ such that every graph $G$ that is both $H$-free and $\overline{H}$-free
satisfies
$\alpha(G)\omega(G)\ge |G|^c$.
\end{thm}

By a theorem of R\"odl~\cite{rodl}, in order to prove \ref{forestsymm} in general, it is enough to prove it for ``sparse'' graphs $G$,
graphs $G$ with maximum degree at most $c|G|$ (for any convenient constant $c$); this argument 
is given in~\cite{cats}. But for sparse graphs, a stronger statement is true (again, a conjecture of~\cite{cats}), that we do 
not need to exclude $\overline{H}$, and the complete
option is no longer needed:
\begin{thm}\label{mainthm}
For every forest $H$ there exists $\epsilon>0$ such that for every $H$-free graph $G$ with $|G|>1$, either
\begin{itemize}
\item some vertex has degree at least $\epsilon |G|$; or
\item there exist disjoint $A,B\subseteq V(G)$ with $|A|,|B|\ge \epsilon |G|$, anticomplete.
\end{itemize}
\end{thm}
(And once again, sparse random graphs show that no non-forest $H$ has this property.)
This is our main result; the derivations of \ref{forestsymm} and \ref{treeEH} from \ref{mainthm} are given in~\cite{cats}.
The proof of \ref{mainthm} is given at the end of section \ref{sec:proofend}.

By working with sparse graphs, we are following a well-trodden path: the papers~\cite{lagoutte, hooks, cats} mentioned above all follow this route, 
by applying R\"odl's theorem. But from this point on,
our method seems to be new.

Let us say a graph $G$ is {\em $\epsilon$-coherent}, where $\epsilon>0$, if:
\begin{itemize}
\item $|G|>1$;
\item every vertex has degree less than $\epsilon |G|$; and
\item there do not exist disjoint anticomplete subsets $A,B\subseteq V(G)$ with $|A|,|B|\ge \epsilon |G|$.
\end{itemize}
(It follows easily that $|G|> \epsilon^{-1}$.)
Thus, \ref{mainthm} is the assertion that for every forest $H$, there exists $\epsilon>0$ such that 
every $\epsilon$-coherent graph contains $H$.

We mention two other papers on $\epsilon$-coherent graphs. First, Bonamy, Bousquet and Thomass\'e~\cite{bonamy} proved that
for every $k$, there exists $\epsilon>0$ such that every $\epsilon$-coherent graph has an induced cycle of length at least $k$;
and second, we proved in~\cite{sparse} the much more general result that for any graph $H$ there exists $\epsilon$ such that every 
$\epsilon$-coherent graph has an induced subgraph that is a subdivision of $H$.

\section{Blockades, and a sketch of the proof}

A {\em blockade} in $G$ means a sequence $(B_i: i\in I)$ of pairwise disjoint nonempty subsets of $V(G)$, where $I$ is a 
finite set of integers.
Its {\em length} is $|I|$,
and the minimum of $|B_i|\;(i\in I)$ is its {\em width}. 
We call the sets $B_i$ {\em blocks} of the blockade. 
(What matters is that the blocks are not too small. We could shrink the larger ones to make them all the same size.)
We are interested in blockades of some fixed length in which each block contains
linearly many vertices of $G$.

Here are two useful ways to make smaller blockades from larger.
First, if $\mathcal{B}=(B_i:i\in I)$ is a blockade, let $I'\subseteq I$; then $(B_i:i\in I')$ is a blockade, of smaller length 
but of at least the same width, and we call it a {\em sub-blockade} of $\mathcal{B}$.
Second, for each $i\in I$ let $B_{i}'\subseteq B_{i}$ be nonempty; then
the sequence $(B_i':i\in I)$ is a blockade, of the same length but possibly of smaller width, and we call it a {\em contraction} 
of $\mathcal{B}$. A contraction of a sub-blockade (or equivalently, a sub-blockade of a contraction) we call a {\em minor} of $\mathcal{B}$.

Let us give an idea of the proof of \ref{mainthm}. If $\mathcal{B}$ is a blockade in $G$, we say that an induced subgraph $H$ of $G$
is {\em $\mathcal{B}$-rainbow} if each vertex of $H$ belongs to some block of $\mathcal{B}$, and no two of them belong to the 
same block. We will prove a stronger theorem:

\begin{thm}\label{rainbow}
For every tree $T$, there exist $d>0$ and an integer $K$, such that, 
for every graph $G$ with a blockade $\mathcal{B}$ of length at least $K$, if $G$ is $\frac{W}{d|G|}$-coherent where $W$ is the width
of $\mathcal{B}$, then
there is a $\mathcal{B}$-rainbow copy of $T$ in $G$.
\end{thm}
(By a {\em copy} of a graph $T$ in $G$ we mean an induced subgraph of $G$ isomorphic to $T$.)
It follows that \ref{rainbow} also holds when $T$ is a forest, because every forest is an induced subgraph of a tree, but for inductive
purposes we state it in terms of trees.
We will prove \ref{rainbow} by induction on $|T|$. Let $v$ be a leaf of $T$, and let $K',d'$ satisfy \ref{rainbow} for the 
forest $T'=T\setminus \{v\}$. Choose $K$ much bigger than $K'$, and $d$ much bigger than $d'$. Let $T'$ have $t$ vertices.

Let $G, \mathcal{B}, w$ be as in the theorem; then we know that for every sub-blockade $\mathcal{B}'$  of $\mathcal{B}$ of length $K'$
there is a $\mathcal{B}'$-rainbow copy of $T'$, which therefore uses some $t$ of the blocks of $\mathcal{B}'$. However, we do not know 
which particular set of $t$ blocks will be used; and even if we did, we still would not know which block contains 
which vertex of $T'$. 

But this can be repaired, very conveniently.
Fix an order of the vertices of $T'$; and for $1\le r_1<\cdots< r_{t}\le K$, say the $t$-tuple $(r_1\ll r_{t})$ is {\em good}
if there is a rainbow copy of $T'$, where for each $i$, the $i$th vertex of $T'$ (in the chosen order) belongs to the block $B_{r_i}$.
Some $t$-tuples may be good and some not; but by Ramsey's theorem, we can find a sub-blockade $\mathcal{B}'$ of $\mathcal{B}$, 
still with 
very large length, such that either every $t$-tuple is good or none of them are. Now we repeat this for all possible orderings of 
the vertices of $T'$. (Note that doing this for a second ordering will not spoil the property we arranged for the first ordering,
as we are dropping to a sub-blockade of $\mathcal{B}'$.)
We can do the same for all ordered trees with at most $t$ vertices.
By this process we produce a sub-blockade of the original blockade, with length much smaller than before (but still as big as we want),
such that if a tree with at most $t$ vertices is present, rainbow, in some subsequence of its blocks, then a copy also appears
rainbow in {\em every} subsequence of the same length, with the vertices of the tree in the same order. This is called being
``support-uniform''. We will not need the 
original blockade $\mathcal{B}$ any more; so for this sketch, to avoid proliferating symbols, let us abuse notation and call the new blockade $\mathcal{B}$, and 
let its length be $K$. Let $\mathcal{B}=(B_i:i\in I)$ where $I=\{1\ll K\}$ say.

Now we have a plentiful supply of copies of $T'$. Indeed, $T'$ must appear $\mathcal{B}$-rainbow somewhere, and therefore it appears 
$\mathcal{B}$-rainbow everywhere and in the same order. Recall that $T'=T\setminus \{v\}$, and let $u$ be the neighbour of $v$ in $T$. We need 
to add back the missing leaf $v$ to one of these many copies of $T'$. How can we do that? 

If $X,Y\subseteq V(G)$ are disjoint, we say that $X$ {\em covers} $Y$ if every vertex in $Y$ has a neighbour in $X$; and $X$ {\em misses}
$Y$ if no vertex in $Y$ has a neighbour in $X$.
Fix an ordering of $V(T')$ such that $T'$ appears in that order, rainbow in every sub-blockade of $\mathcal{B}$ of length $t$; 
and let $u$ be the $k$th vertex in this ordering. Thus there are $k-1$ vertices of $T'$ that are earlier than $u$, and $t-k$ 
that are later.
Suppose that for some $j$ with $1\le j\le K$, there is a subset $X$ of $V(G)$ that 
\begin{itemize}
\item covers at least a positive fraction 
(say $1/100$) of $B_j$; 
\item misses at least $1/100$ of $B_h$ for $k-1$ values of $h<i$ different from $j$, say $r_1<\cdots < r_{k-1}$; 
\item misses at least $1/100$ of $B_h$ for $t-k$ values of $h>i$ different from $j$, say $r_{k+1}<\cdots < r_{t}$; and
\item is included in the union of the blocks $B_h$ where $h\ne r_1\ll  r_{k-1}, j, r_{k+1}\ll r_{t}$.
\end{itemize}
We know that there is a $\mathcal{B}$-rainbow copy of $T'$ that appears in the right order in the blocks 
$$B_{r_1}, B_{r_2}\ll B_{r_{k-1}}, B_j, B_{r_{k+1}}\ll B_{r_{t}},$$
and in particular vertex $u$ appears in block $B_j$. If we could arrange that $X$ covers all of $B_j$, and missed all of the other blocks
$B_{r_1}, B_{r_2}\ll B_{r_{k-1}}, B_{r_{k+1}}\ll B_{r_{t}}$, then we can recover the missing leaf (use a vertex in $X$ adjacent
to $u$). But this is not the case: all we know is that $X$ covers $1/100$ of $B_j$, and so on. To fix this, we need to know that if we shrink
each block of $\mathcal{B}$ to $1/100$ of its present size, then rainbow copies of this ordering of $T'$ are still present everywhere.

But that is another thing we could arrange. If by shrinking the blocks to $1/100$ of their present size, we can stop some $t$-tuple 
being good, then do so, and repeat for all $t$-tuples and all orderings of $T'$. Eliminate as many orderings of $T'$
from as many $t$-tuples of blocks as possible. When this stops (and it will stop after a constant number of steps, because
$K$ is only a constant), we have spoiled the support-uniformity property; 
but we can repair it by doing the Ramsey argument again. (Actually, we just do it once, but after the shrinking process, not before.)

So now we have constructed a very interesting blockade (let us call it $\mathcal{B}$ again): its length is still a big constant; its width ($w$ say) 
is not too small a fraction of $|G|$; every ordered 
tree with $t$ vertices that appears rainbow in some sub-blockade of length $t$ appears in every contraction of this blockade 
whose width is at least $w/100$ (we call this ``support-invariance''); and appears in every sub-blockade of length $t$; and
we can assume that there is no $j,X$ satisfying the bullets above (because otherwise we have a rainbow copy of $T$).

Let us partition the blocks of $\mathcal{B}$ into many intervals, each including many blocks; we get another blockade $\mathcal{B}'$ say, still support-uniform and support-invariant, and if we can
find a $\mathcal{B}'$-rainbow copy of $T$, then it is also $\mathcal{B}$-rainbow.
And the fact that no pair $j,X$ satisfies the three bullets above in $\mathcal{B}$ translates into something nicer for $\mathcal{B}'$,
an important property we call ``concavity'': there do not exist distinct blocks $B_h', B_i', B_j'$ with $h<i<j$, such that 
for some $X$ included in the union of the other blocks, $X$ covers at least $1/50$ of $B_i'$ and misses at least $1/50$ of $B_h', B_j'$. 
(For this step to work, we need first to arrange that all the blocks of $\mathcal{B}'$ have the same size, but we can do that.)

That seems to be the limit of what we can get by this method, using the inductive hypothesis on $V(T)$ to obtain a well-positioned 
copy of $T'$,
in a position where the missing vertex can be replaced. For the remainder of the proof we need a different approach.
So far we have not used the hypothesis 
about coherence, but now that will come into play. In the next section, we will show that for any tree $T$, if a graph $G$ has a support-uniform, 
support-invariant concave blockade $\mathcal{B}$ 
of sufficient length, and 
width at least some $w$, and $G$ is $w/(d|G|)$-coherent where $d$ is some constant independent of $G$, then there is a 
$\mathcal{B}$-rainbow copy of $T$. 
Then in section \ref{sec:proofend}, we use that result, combined with the argument sketched in this section, to prove \ref{rainbow} and hence \ref{mainthm}.

\section{Using concavity}

An {\em ordered graph} is a graph $J$ together with a linear order of its vertex set. Isomorphism of
ordered graphs is defined in the natural way.
Let $\mathcal{B}=(B_i:i\in I)$ be a blockade in $G$.
We recall that an induced subgraph $H$ of $G$ is {\em $\mathcal{B}$-rainbow}
if each vertex of $H$ belongs to some block of $\mathcal{B}$, and no two vertices belong to the same block.
If $H$ is a $\mathcal{B}$-rainbow induced subgraph of $G$, there is an associated 
ordering $<$ of the vertex set of $H$, where we say 
$u<v$ if $u\in B_i$ and $v\in B_j$
with $i<j$. We call the ordered graph given by this ordering the {\em $\mathcal{B}$-ordering} of $H$.
If $J$ is an ordered graph, and the $\mathcal{B}$-ordering of $H$ is isomorphic to $J$, we say that $H$ is a {\em copy} of $J$.

If $H$ is a $\mathcal{B}$-rainbow induced subgraph,
its {\em support} is the set of all $i\in I$ such that $V(H)\cap B_i\ne \emptyset$.
If $J$ is an ordered graph, we define the {\em trace} of $J$ (relative to $\mathcal{B}$) to be the set of supports of all $\mathcal{B}$-rainbow copies of $J$.
If $\tau\ge 1$ an integer, we say $\mathcal{B}$ is {\em $\tau$-support-uniform} if
for every ordered tree $J$ with $|J|\le \tau$, either the 
trace of $J$ is empty, or it consists of all subsets of $I$ of cardinality $|J|$.

Let $0< \kappa\le 1$ and $\tau\ge 1$.
We say $\mathcal{B}$ is
{\em $(\kappa,\tau)$-support-invariant} if 
for every contraction $\mathcal{B}'=(B_i':i\in I)$ of $\mathcal{B}$ of width at least $\kappa$ times the width of $\mathcal{B}$,
and for every ordered tree $J$ with $|J|\le \tau$,
the trace of $J$ relative to $\mathcal{B}$ equals the trace of $J$ relative to $\mathcal{B'}$.

Say a blockade is {\em equicardinal} if all its blocks have the same cardinality.
Let $\mathcal{B}=(B_1\ll B_k)$ be an equicardinal blockade of width $w$, and let $0\le \lambda\le  1/2$. For $1\le i\le k$, a subset $X$ of
$V(G)\setminus B_i$ is said to {\em $\lambda$-cover}
$B_i$ if there are at least $\lambda w$ vertices in $B_i$ with a neighbour in $X$, and to {\em $\lambda$-miss}
$B_i$ if there are at least $\lambda w$ vertices in $B_i$ with no neighbour in $X$.
Since $\lambda\le 1/2$ and $|B_i|\ge w$, $X$ either $\lambda$-covers $B_i$
or $\lambda$-misses $B_i$, but it might do both.

We denote the union of all the blocks of a blockade $\mathcal{B}$ by $V(\mathcal{B})$.
An equicardinal blockade $\mathcal{B}=(B_i:i\in I)$ is {\em $\lambda$-concave} if it has the following very strong property:
for all $h_1,h_2,h_3\in I$ with $h_1<h_2<h_3$, there is no subset 
$X\subseteq V(\mathcal{B})\setminus (B_{h_1}\cup B_{h_2}\cup B_{h_3})$ that
$\lambda$-covers $B_{h_2}$ and $\lambda$-misses $B_{h_1}$ and $B_{h_3}$.

A {\rm rooted graph} $H$ is a pair $(H^-,r(H))$, where $H^-$ is a graph and $r(H)\in V(H^-)$; we call $r(H)$ the {\em root}.
If $H_1,H_2$ are rooted graphs, by an isomorphism between them we mean an isomorphism
between $H_1^-$ and $H_2^-$ that takes root to root.
%Define rooted subtree
If $\delta\ge 2$ is an integer, let $T(\delta,0)$ be the rooted tree with one vertex
(thus, $\delta$ is irrelevant, but this will be convenient).
If $\eta\ge 1$ and $\delta\ge 2$ are integers, we denote by $T(\delta, \eta)$ the rooted tree with the properties that
\begin{itemize}
\item every vertex has degree $\delta+1$ or $1$, except the root, which has degree $\delta$; and
\item for each vertex of degree one, its distance from the root is exactly $\eta$.
\end{itemize}
Thus for $\eta\ge 1$, $T(\delta, \eta)$  is formed by taking the disjoint union of $\delta$ copies of $T(\delta, \eta-1)$,
and adding a new vertex adjacent to all the roots, and making this vertex the new root.
Clearly every tree is isomorphic to
a subtree of $T(\delta,\eta)$ for sufficiently large $\delta$ and
$\eta$.
\begin{figure}[H]
\centering

\begin{tikzpicture}[scale=0.8,auto=left]
\tikzstyle{every node}=[inner sep=1.5pt, fill=black,circle,draw]
\node (a) at (0,0) {};
\node (b1) at (-2,-1) {};
\node (b2) at (2,-1) {};
\node (c2) at (-1,-2) {};
\node (c1) at (-3,-2) {};
\node (c3) at (1,-2) {};
\node (c4) at (3,-2) {};
\node (d1) at (-3.5,-3) {};
\node (d2) at (-2.5,-3) {};
\node (d3) at (-1.5,-3) {};
\node (d4) at (-.5,-3) {};
\node (d5) at (.5,-3) {};
\node (d6) at (1.5,-3) {};
\node (d7) at (3.5,-3) {};
\node (d8) at (2.5,-3) {};

\foreach \from/\to in {a/b1,a/b2,b1/c1,b1/c2,b2/c3,b2/c4,c1/d1,c1/d2,c2/d3,c2/d4,c3/d5,c3/d6,c4/d7,c4/d8}
\draw [-] (\from) -- (\to);

\tikzstyle{every node}=[]
\draw[above] (a) node []           {root};

\end{tikzpicture}
\caption{$T(2,3)$.} \label{fig:T(2,3)}
\end{figure}

Let $\mathcal{B}=(B_1\ll B_K)$ be a blockade in $G$. 
A $\mathcal{B}$-rainbow induced rooted subgraph $H$ of $G$ is {\em $\mathcal{B}$-left-rainbow} 
if
$j\ge i$ for all $j\in \{1\ll K\}$ with $V(H)\cap B_j\ne \emptyset$, where the root of $H$ belongs to $B_i$.
We define {\em $\mathcal{B}$-right-rainbow} similarly, requiring $j\le i$ instead.

\begin{thm}\label{usemonotone}
Let $\delta\ge 2$ and $\eta\ge 0$ be integers. Let 
$0<\lambda \le 2^{-9\delta}\delta^{-1-\eta}$, let $\tau= \delta^{1+\eta}$, and let $\epsilon>0$.
Let $G$ be an $\epsilon$-coherent graph with an equicardinal blockade
$\mathcal{B}$ of length at least $6\delta^{\eta+2}$ and width at least $2^{9\delta}\epsilon|G|$, such that
$\mathcal{B}$ is $\lambda$-concave, $\tau$-support-uniform and $(2^{-9\delta},\tau)$-support-invariant.
Then $G$ contains $T(\delta,\eta)$.
\end{thm}
\Proof
Let $\mathcal{B}=(B_1\ll B_K)$, and let $W$ be its width.
Thus $|B_i|=W$ for $1\le i\le K$.
Choose $\alpha\ge 0$ maximum such that there is a $\mathcal{B}$-left-rainbow copy of $T(\delta,\alpha)$,
and define $\beta$ similarly for $\mathcal{B}$-right-rainbow. We suppose for a contradiction that there is no $\mathcal{B}$-rainbow copy of $T(\delta,\eta)$, and so
$\alpha,\beta<\eta$; and by reversing the blockade if necessary we may assume that $\alpha\le \beta$.
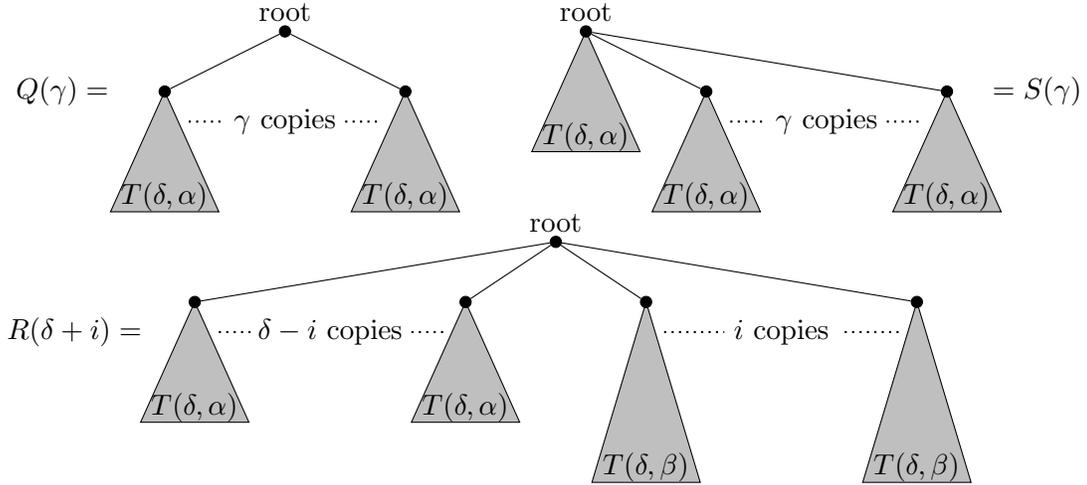
\begin{figure}[h!]
\centering

\begin{tikzpicture}[scale=0.8,auto=left]
\draw[fill= lightgray] (-2,-1)--(-2.9,-3)--(-1.1,-3)--(-2,-1);
\draw[fill= lightgray] (2,-1)--(2.9,-3)--(1.1,-3)--(2,-1);
\tikzstyle{every node}=[inner sep=1.5pt, fill=black,circle,draw]
\node (a) at (0,0) {};
\node (b) at (-2,-1) {};
\node (c) at (2,-1) {};
\tikzstyle{every node}=[]
\draw[above] (a) node []           {root};
\node at (-2,-2.75) {$T(\delta,\alpha)$};
\node at (2,-2.75) {$T(\delta,\alpha)$};
\node at (0,-1.5) {$\gamma$ copies};

\draw[-] (a)--(b);
\draw[-] (a)--(c);

\draw[dotted,thick] (-1.6,-1.5)--(-1,-1.5);
\draw[dotted,thick] (1,-1.5)--(1.6,-1.5);
\node at (-3.7,-1) {$Q(\gamma)=$};

\begin{scope}[shift ={(9,0)}]
\draw[fill= lightgray] (-4,0)--(-4.9,-2)--(-3.1,-2)--(-4,0);
\draw[fill= lightgray] (-2,-1)--(-2.9,-3)--(-1.1,-3)--(-2,-1);
\draw[fill= lightgray] (2,-1)--(2.9,-3)--(1.1,-3)--(2,-1);
\tikzstyle{every node}=[inner sep=1.5pt, fill=black,circle,draw]
\node (a) at (-4,0) {};
\node (b) at (-2,-1) {};
\node (c) at (2,-1) {};
\tikzstyle{every node}=[]
\draw[above] (a) node []           {root};
\node at (-2,-2.75) {$T(\delta,\alpha)$};
\node at (2,-2.75) {$T(\delta,\alpha)$};
\node at (-4,-1.75) {$T(\delta,\alpha)$};
\node at (0,-1.5) {$\gamma$ copies};

\draw[-] (a)--(b);
\draw[-] (a)--(c);

\draw[dotted,thick] (-1.6,-1.5)--(-1,-1.5);
\draw[dotted,thick] (1,-1.5)--(1.6,-1.5);

\node at (3.5,-1) {$=S(\gamma)$};
\end{scope}

\begin{scope}[shift ={(2,-3.5)}]
\draw[fill= lightgray] (-3.5,-1)--(-4.4,-3)--(-2.6,-3)--(-3.5,-1);
\draw[fill= lightgray] (1,-1)--(1.9,-3)--(.1,-3)--(1,-1);
\tikzstyle{every node}=[inner sep=1.5pt, fill=black,circle,draw]
\node (b) at (-3.5,-1) {};
\node (c) at (1,-1) {};
\tikzstyle{every node}=[]
\node at (-3.5,-2.75) {$T(\delta,\alpha)$};
\node at (1,-2.75) {$T(\delta,\alpha)$};
\node at (-1.25,-1.5) {$\delta-i$ copies};

\draw[dotted,thick] (-3.1,-1.5)--(-2.5,-1.5);
\draw[dotted,thick] (0.1,-1.5)--(0.7,-1.5);
\node at (-5.5,-1.5) {$R(\delta+i)=$};
\end{scope}

\begin{scope}[shift ={(8,-3.5)}]
\draw[fill= lightgray] (-2,-1)--(-2.9,-4)--(-1.1,-4)--(-2,-1);
\draw[fill= lightgray] (2.5,-1)--(3.4,-4)--(1.6,-4)--(2.5,-1);
\tikzstyle{every node}=[inner sep=1.5pt, fill=black,circle,draw]
\node (a) at (-3.5,0) {};
\node (d) at (-2,-1) {};
\node (e) at (2.5,-1) {};
\tikzstyle{every node}=[]
\draw[above] (a) node []           {root};
\node at (-2,-3.75) {$T(\delta,\beta)$};
\node at (2.5,-3.75) {$T(\delta,\beta)$};
\node at (.25,-1.5) {$i$ copies};

\draw[-] (a)--(d);
\draw[-] (a)--(e);

\draw[dotted,thick] (-1.7,-1.5)--(-.7,-1.5);
\draw[dotted,thick] (1.3,-1.5)--(2.2,-1.5);

\end{scope}
\draw[-] (a)--(b);
\draw[-] (a)--(c);
\end{tikzpicture}
\caption{$Q(\gamma)$, $R(\gamma)$ and $S(\gamma)$.} \label{R(gamma)}
\end{figure}

We need three special rooted trees:
\begin{itemize}
\item For $0\le \gamma\le \delta$, let $Q(\gamma)$ be obtained from the disjoint
union of
$\gamma$ copies of $T(\delta,\alpha)$ by adding a new root adjacent to the old roots.
\item  For $0\le \gamma\le \delta$ let $R(\gamma)=Q(\gamma)$, and for $0<i\le \delta$ let $R(\delta+i)$ be obtained from the disjoint union of
$\delta-i$ copies of $T(\delta,\alpha)$ and $i$ copies of $T(\delta,\beta)$ by adding a new root adjacent to all
the old roots. (Thus, as we go from $R(0)$ to $R(2\delta)$ we first add copies of $T(\delta,\alpha)$ one at a time, and then replace them
with $T(\delta,\beta)$ one at a time, finishing with $T(\delta,\beta+1)$.)
\item For $0\le \gamma\le \delta$, let $S(\gamma)$ be obtained from the disjoint union of $\gamma+1$ copies of
$T(\delta,\alpha)$ by making the root $v$ of the first copy adjacent
to all other roots, and making $v$ the new root.
\end{itemize}

Let $\mathcal{B}'=(B_i':i\in I)$ be a minor of $\mathcal{B}$, with $B_i'\subseteq B_i$ for each $i\in I$.
We say $\mathcal{B}'$ is {\em $(\gamma_1,\gamma_2)$-anchored} in $\mathcal{B}$ if
there exist $k, Y$ such that:
\begin{itemize}
\item $k$ is an integer with $1\le k< K$ such that $I=\{1\ll k\}\cup \{K\}$;
\item $Y$ is a subset of $\bigcup_{k< j< K} B_j$;
\item $Y$ is anticomplete to $B_{i}'$ for all $i\in I\setminus \{1,K\}$;
\item for every $v\in B_{1}'$ there is a $\mathcal{B}$-left-rainbow copy of $Q(\gamma_1)$ in $G[Y\cup \{v\}]$ with root $v$;
and
\item for every $v\in B_{K}'$ there is a $\mathcal{B}$-right-rainbow copy of $R(\gamma_2)$ in $G[Y\cup \{v\}]$ with root $v$.
\end{itemize}
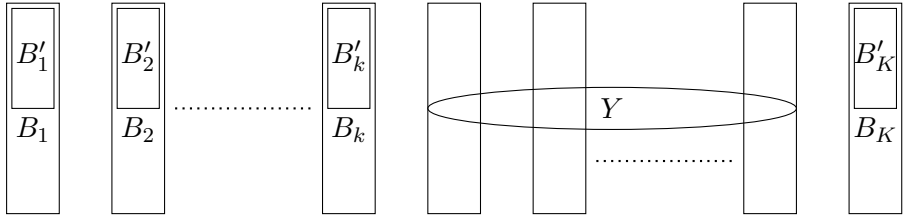
\begin{figure}[h!]
\centering

\begin{tikzpicture}[scale=0.7,auto=left]
\tikzstyle{every node}=[]

\draw (0,0) rectangle (1,4);
\draw (2,0) rectangle (3,4);
\draw (6,0) rectangle (7,4);
\draw (8,0) rectangle (9,4);
\draw (10,0) rectangle (11,4);
\draw (14,0) rectangle (15,4);
\draw (16,0) rectangle (17,4);

\draw (0.1,2) rectangle (.9,3.9);
\draw (2.1,2) rectangle (2.9,3.9);
\draw (6.1,2) rectangle (6.9,3.9);
\draw (16.1,2) rectangle (16.9,3.9);

\draw[dotted, thick] (3.2,2) -- (5.8,2);
\draw[dotted, thick] (11.2,1) -- (13.8,1);

\draw (11.5,2) ellipse (3.5 and .4);

\node at (.5,1.6) {$B_1$};
\node at (2.5,1.6) {$B_2$};
\node at (6.5,1.6) {$B_{k}$};
\node at (16.5,1.6) {$B_K$};

\node at (.5,3) {$B_1'$};
\node at (2.5,3) {$B_2'$};
\node at (6.5,3) {$B_k'$};
\node at (16.5,3) {$B_K'$};
\node at (11.5,2) {$Y$};
\end{tikzpicture}
\caption{Figure for ``anchored''.} \label{fig:anchored}
\end{figure}
Anchored blockades will be important. The minor
$\mathcal{B}'$ will inherit the concavity of $\mathcal{B}$; and since $Y$ only sees the
two blocks at either end of $\mathcal{B}'$, we can grow trees rooted in $B_1'$ or $B_K'$ and
then extend them using $Y$.  This will give us most of what we want.  However,
there will be one awkward case that will make the end of the argument a
little more complicated (and is one reason for having three types of special
tree).

Choose $\gamma_0\ge 0$ maximum such that there exist $\gamma_1,\gamma_2\ge 0$ with $\gamma_1+\gamma_2=\gamma_0$
and  a $(\gamma_1, \gamma_2)$-anchored minor $\mathcal{B}'$
of $\mathcal{B}$ of length at least $K-2\delta^{\eta+1}\gamma_0$ and
width at least $W2^{-3\gamma_0}$. (This is possible, because for $\gamma_0=0$ we can take $\mathcal{B}'=\mathcal{B}$ and $k=K-1$.)
Let $W'\ge W2^{-3\gamma_0}$ be its width. We may assume that $\mathcal{B}'$ is equicardinal (by
shrinking any blocks that are larger: all the defining properties of $\mathcal{B}'$ are preserved).
Choose $\gamma_3$ maximum such that there is a $\mathcal{B}$-left-rainbow copy of $S(\gamma_3)$ (this is possible since
there is a $\mathcal{B}$-left-rainbow copy of $S(0)=T(\alpha,\delta)$). We claim that:
\\
\\
(1) {\em $\gamma_1,\gamma_3\le \delta-1$, and $\gamma_2\le 2\delta-1$. Consequently $\gamma_0\le 3\delta-2$, and
so $W'\ge W2^{-9\delta+6}\ge 64\epsilon|G|$.}
\\
\\
Since $Q(\delta)$ is isomorphic to
$T(\delta, \alpha+1)$, it follows from the choice of $\alpha$ that there is no $\mathcal{B}$-left-rainbow
copy of $Q(\delta)$. On the other hand there is a $\mathcal{B}$-left-rainbow copy of $Q(\gamma_1)$,
from the definition of ``$(\gamma_1, \gamma_2)$-anchored''; so
$\gamma_1<\delta$. Similarly $\gamma_2< 2\delta$, since $R(2\delta)$ is isomorphic to $T(\delta, \beta+1)$.
Also $\gamma_3<\delta$ from the maximality of $\alpha$, since $S(\delta)$ contains $T(\delta,\alpha+1)$.
This proves (1).

\bigskip

Now $\mathcal{B}'$ is $(\gamma_1, \gamma_2)$-anchored; let $Y, k$ be as in the definition of ``anchored'', and let
$\mathcal{B}' = (B_{i}':i\in \{1\ll k\}\cup \{K\})$, where $B_i'\subseteq B_i$ for each $i$.
Let $S(\gamma_3)$ have $s$ vertices, and let $T(\delta,\beta)$
have $t$ vertices. Define $h=k+1-s-t$.  Next we show that:
\\
\\
(2) {\em $h\ge K-2\delta^{\eta+1}(\gamma_0+1)\ge 1$.}
\\
\\
Since $\mathcal{B}'$ has length $k+1$, it follows that $k+1\ge K-2\delta^{\eta+1}\gamma_0$.
Hence $h=k+1-s-t\ge K-2\delta^{\eta+1}(\gamma_0+1)$ since $s,t\le \delta^{\eta+1}$.
But
$$K\ge 6\delta^{\eta+2}>2\delta^{\eta+1}(3\delta-1)\ge 2\delta^{\eta+1}(\gamma_0+1),$$
by (1), and 
it follows that $h>0$.
This proves (2).
%K>\delta{\eta+1}(3\delta)+1
%K\ge 6\delta^{\eta+2}

\bigskip

Let $r=\lceil W'-2^{-9\delta} W\rceil$. By (1), $W'\ge W2^{-9\delta+6}$, and so $r\ge 63\cdot 2^{-9\delta} W$.
Since $W\ge 2^{9\delta}\epsilon|G|$, this implies that $r\ge 63 \epsilon |G|$. We show next that:
%W_0-\mu W\ge 0
\\
\\
(3) {\em There are $r$ copies $E_1\ll E_r$ of $S(\gamma_3)$,
pairwise vertex-disjoint and each $(B_i':h\le i\le h+s-1)$-left-rainbow; and there are $r$ copies $F_1\ll F_r$ of $T(\delta,\beta)$,
pairwise vertex-disjoint and $(B_i':h+s\le i\le  k)$-right-rainbow.}
\\
\\
Since there is a $\mathcal{B}$-left-rainbow copy of $S(\gamma_3)$, and $\mathcal{B}$ is $\tau$-support-uniform,
and $|S(\gamma_3)|\le \tau$,
there is such a copy that is $(B_i:h\le i\le h+s-1)$-left-rainbow.
Choose $r'\le r$ maximum such that there are $r'$ pairwise disjoint copies of $S(\gamma_3)$,
pairwise vertex-disjoint and $(B_i':h\le i\le h+s-1)$-left-rainbow. By
removing the vertices of these copies from the blocks $B_i'\; (i\in \{h\ll h+s-1\})$,
we obtain a contraction of $(B_i:h\le i\le h+s-1)$ of width $W'-r'$ in which there is no left-rainbow copy of $S(\gamma_3)$.
But $(B_i:h\le i\le h+s-1)$ is
$(2^{-9\delta},\tau)$-support-invariant,
and so $W'-r'<2^{-9\delta} W$, that is, $r'=r$. This proves the first assertion, and the second follows similarly. This proves (3).

\bigskip
Thus the collection of rooted trees $E_i$ and
$F_i$ form a large ``matching'' in the interval $(B_h'\ll B_k')$, with roots in the
blocks $B_h'$ and $B_k'$ that lie at the left and right hand ends.
\begin{figure}[H]
\centering

\begin{tikzpicture}[scale=0.7,auto=left]
\tikzstyle{every node}=[]

\draw (1,0) rectangle (2,4);
\draw (1.1,1) rectangle (1.9,3.9);
\node at (1.5,.6) {$B_1$};
\node at (1.5,2.5) {$B_1'$};

\draw (4,0) rectangle (5,4);
\draw (4.1,1) rectangle (4.9,3.9);

\draw (6,0) rectangle (7,4);
\draw (6.1,1) rectangle (6.9,3.9);
\node at (6.5,.6) {$B_{h}$};
\node at (6.5,1.5) {$B_h'$};

\draw (10,0) rectangle (11,4);
\draw (10.1,1) rectangle (10.9,3.9);

\draw (12,0) rectangle (13,4);
\draw (12.1,1) rectangle (12.9,3.9);

\draw (16,0) rectangle (17,4);
\draw (16.1,1) rectangle (16.9,3.9);
\node at (16.5,.6) {$B_{k}$};
\node at (16.5,1.5) {$B_k'$};

\draw (18,0) rectangle (19,4);

\draw (21,0) rectangle (22,4);

\draw (23,0) rectangle (24,4);
\draw (23.1,1) rectangle (23.9,3.9);
\node at (23.5,.6) {$B_K$};
\node at (23.5,2.5) {$B_K'$};

\draw[dotted, thick] (2.2,1.2) -- (3.8,1.2);
\draw[dotted, thick] (7.2,1.2) -- (9.8,1.2);
\draw[dotted, thick] (13.2,1.2) -- (15.8,1.2);
\draw[dotted, thick] (19.2,1.2) -- (20.8,1.2);

\draw (20,2) ellipse (2 and .3);
\node at (20,2) {$Y$};

\draw (8.5,3.7) ellipse (2.4 and .3);
\node at (8.5,3.7) {$E_1$};
\draw (8.5,2.2) ellipse (2.4 and .3);
\node at (8.5,2.2) {$E_r$};
\draw (14.5,3.7) ellipse (2.4 and .3);
\node at (14.5,3.7) {$F_1$};
\draw (14.5,2.2) ellipse (2.4 and .3);
\node at (14.5,2.2) {$F_r$};

\draw[dotted, thick] (8.5,3.2) -- (8.5,2.7);
\draw[dotted, thick] (14.5,3.2) -- (14.5,2.7);

\tikzstyle{every node}=[inner sep=1.5pt, fill=black,circle,draw]
\node (a) at (6.5,3.7) {};
\node (b) at (6.5,2.2) {};
\node (c) at (16.5,3.7) {};
\node (d) at (16.5,2.2) {};

\end{tikzpicture}
\caption{The position of the rooted trees $E_i$ and $F_i$.} \label{fig:anchored2}
\end{figure}
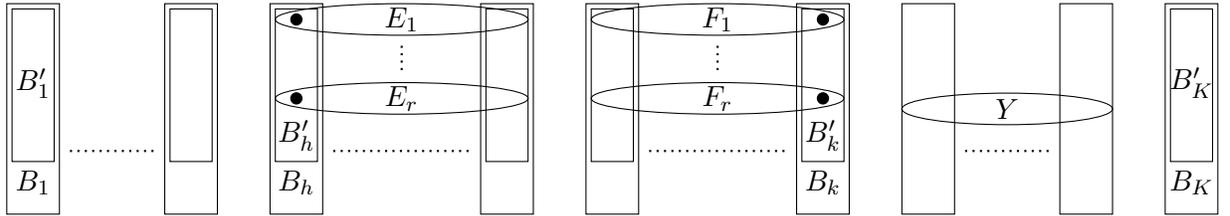

\bigskip

In order to grow a
tree, we will want to join the roots of $E_i$ or $F_i$ to the blocks $B_1'$ and
$B_K'$; however we need to consider possible adjacencies to other vertices in
$E_i$ and $F_i$.

For $v\in B_1'\cup B_K'$, and $1\le i\le r$, we say
\begin{itemize}
\item  $v$ {\em meets} $E_i\cup F_i$ if $v$ is adjacent to some vertex of
$E_i\cup F_i$;
\item $v$ meets $E_i\cup F_i$ {\em internally} if $v$ is adjacent to some vertex of $E_i\cup F_i$
that is not the root of $E_i$ or $F_i$ (and possibly $v$ is also adjacent to one or both roots);
\item $v$ meets $E_i\cup F_i$ {\em properly} if $v$ is adjacent to one or both of the roots of $E_i, F_i$, but to no other
vertices of $E_i\cup F_i$, that is, if $v$ meets $E_i\cup F_i$ and does not meet $E_i\cup F_i$ internally.
\end{itemize}

For $X\subseteq B_1'\cup B_K'$, let $a(X)$ be the number of $i\in \{1\ll r\}$ such that some vertex in $X$ meets $E_i\cup F_i$,
and let $b(X)$ be the number of $i\in \{1\ll r\}$ such that some vertex in $X$ meets $E_i\cup F_i$ internally. Choose $D\subseteq B_1'\cup B_K'$
maximal such that $a(D)\le r/2$ and $b(D)\ge a(D)/4$. We claim that:
\\
\\
(4) {\em $|D|<\epsilon n$, and $a(D)\le r/2-\epsilon n$.}
\\
\\
There are at least $r/2$ vertices in $B_{h}$
with no neighbour in $D$ (the roots of the trees $E_i$ such that no vertex in $D$ meets $E_i\cup F_i$). Since
$r/2\ge \epsilon n$ and $G$ is $\epsilon$-coherent, it follows that
$|D|< \epsilon n$; and since $r/2\ge \lambda Wn$,
%$r/2\ge \lambda Wn$,
it follows that $D$ $\lambda$-misses $B_{h}$.
Similarly $D$ $\lambda$-misses $B_{h+s+t-1}$, and since
$\mathcal{B}$ is $\lambda$-concave, $D$ does not $\lambda$-cover any of the sets $B_{h+1}\ll B_{h+s+t-2}$. Hence
there are at most $\lambda (s+t-2)Wn$ vertices in $B_{h+1}\cup\cdots \cup B_{h+s+t-2}$ that have neighbours in $D$.
Since
there are at least $b(D)$ such vertices in total, it follows that $\lambda (s+t)Wn\ge b(D)\ge a(D)/4$, and so
$$a(D)\le 4\lambda(s+t)Wn\le 4(2^{-9\delta}\delta^{-1-\eta})(2\delta^{\eta+1})(2^{9\delta}r/63)=8r/63,$$
since $\lambda \le 2^{-9\delta}\delta^{-1-\eta}$, and $s+t\le 2 \delta^{\eta+1}$, and $r\ge 63\cdot 2^{-9\delta} Wn$.
Since $r\ge 63 \epsilon n$, it follows that $8r/63\le r/2-\epsilon n$. This proves (4).
%$r/2\ge \epsilon n$,
%\lceil (W(g)-\mu W)|G|\rceil$\ge \lambda^{1/2}(s+t-2)|G| + \epsilon |G|$??

\bigskip

Let $Z$ be the set of vertices in $(B_1'\cup B_K')\setminus D$ that meet one of $E_1\cup F_1\ll E_r\cup F_r$.
We show that:
\\
\\
(5) {\em $|Z|\ge 2(W'-\epsilon)n$.}
\\
\\
Since $r\ge \epsilon n$,
%$r\ge \epsilon |G|$
and $G$ is $\epsilon$-coherent, there are fewer than $\epsilon n$ vertices in $B_1'\cup B_K'$ that have no neighbour in any of
$E_1\cup F_1\ll E_r\cup F_r$. All the other vertices in $B_1'\cup B_K'$ belong to either $D$ or $Z$, so
$|D|+|Z|+\epsilon n \ge |B_1'\cup B_K'|= 2W'n$.
From (4), this proves (5).

\bigskip
Let $C$ be the set of all $i\in \{1\ll r\}$ such that $D$ is anticomplete to $V(E_i\cup F_i)$. Thus $|C|=r-a(D)$. We claim:
\\
\\
(6) {\em For each $v\in Z$, the number of $i\in C$ such that $v$ meets $E_i\cup F_i$ internally is less than a quarter of
the number of $i\in C$ such that $v$ meets $E_i\cup F_i$.}
\\
\\
Since $a(D)\le r/2-\epsilon n$ by (4), it follows that
$a(D\cup \{v\})\le r/2$, and the maximality of $D$ implies that
$b(D\cup\{v\})< a(D\cup \{v\})/4$. Since $b(D)\ge a(D)/4$, it follows that
$$b(D\cup\{v\})-b(D)<(a(D\cup \{v\})-a(D))/4.$$
But $b(D\cup\{v\})-b(D)$ is at least the number of $i\in C$ such that $v$ meets $E_i\cup F_i$ internally; and $a(D\cup \{v\})-a(D)$
equals the number of $i\in C$ such that $v$ meets $E_i\cup F_i$. 
This proves (6).
%$\lambda^{1/2}(s+t-2)Wn\le r/2-\epsilon n$
\\
\\
(7) {\em We may assume that there is a subset $X\subseteq Z$ with $|X\cap B_1'|\ge |Z\cap B_1'|/2$ and $|X\cap B_K'|\ge |Z\cap B_K'|/2$,
such that for each $v\in X$,
if $i\in C$ is minimum such that $v$ meets $E_i\cup F_i$, then $v$ meets $E_i\cup F_i$ properly.}
\\
\\
Let $v\in Z$, and take a linear order of $C$;
and let $i\in C$ be the first member of $C$ (under this order) such that $v$ meets $E_i\cup F_i$ (there is such a member $i$
from the definition of $Z$). We say $v$
is {\em happy} (under this order), if $v$ meets $E_i\cup F_i$ properly. If we choose the linear order uniformly at random,
the probability that $v$ is happy is more than $3/4$, by (6); and so the expected number of vertices in $Z\cap B_1'$ that are happy
is more than $3|D\cap B_1'|/4$. Hence the probability that at least $|D\cap B_1'|/2$ vertices in $Z\cap B_1'$ are happy is more
than $1/2$, and the same for $Z\cap B_K'$; and so there is a positive probability that both events occur. Hence
there is a linear order of $C$ such that at least
$|Z\cap B_1'|/2$ vertices in $Z\cap B_1'$ are happy, and at least $|Z\cap B_K'|/2$ vertices in $Z\cap B_K'$ are happy.
By renumbering, we may assume that this order
is the natural order of $C$ as a set of integers. This proves (7).

\bigskip

For each $v\in X$, we call the value of $i$ in (7) the {\em happiness} of $v$.
Let $v\in X$ and let $i$ be its happiness. Since $v$ meets $E_i\cup F_i$ properly, it is adjacent
to one or both of the roots of $E_i, F_i$, and has no other neighbours in $V(E_i\cup F_i)$. Also, $v$ belongs to one of $B_1',B_K'$.
Let us say $v$ has
\begin{itemize}
\item type $(1,E)$ if $v\in B_1'$ and $v$ is adjacent to the root of $E_i$;
\item type $(1,F)$ if $v\in B_1'$ and $v$ is adjacent to the root of $F_i$;
\item type $(K,E)$ if $v\in B_K'$ and $v$ is adjacent to the root of $E_i$; and
\item type $(K,F)$ if $v\in B_K'$ and $v$ is adjacent to the root of $F_i$.
\end{itemize}
Every vertex in $X$ has one of these four types (some may have two types).
Now $|Z|\ge 2(W'-\epsilon)n$ by (5), and $|X|\ge |Z|/2$ by (7), and so
$|X|\ge (W'-\epsilon) n\ge 63W'n/64$ by (1).
Thus one of $B_1'\cap X,B_K'\cap X$ has cardinality at least $63W'n/128\ge W'n/4$, so we may
choose $m\le |C|$ minimum such that one of $B_1'\cap X,B_K'\cap X$ contains at least $W'n/4$ vertices with happiness at most $m$.
Consequently there is a set $U\subseteq X$, such that $|U|\ge W'n/8$, and all vertices in $U$ have happiness at most $m$,
and they all have the same type (which,
from now on, we call the ``type of $U$'').
Let
$$Y'=V(E_1)\cup\cdots\cup V(E_m)\cup V(F_1)\cup\cdots\cup V(F_m).$$
We claim that:
\\
\\
(8) {\em For each $i\in \{2\ll h-1\}$, there is a subset of $B_i'$ anticomplete to $Y'$,
of cardinality at least $W'n/8$.}
\\
\\
Let $j\in \{h\ll h+s+t-1\}$.
From the choice of $m$, fewer than $W'n/4$ vertices in $B_1'\cap X$ have happiness less
than $m$; and so at most $W'n/4 +2\epsilon n$
have happiness at most $m$, since those with happiness exactly $m$ are adjacent to one of the roots of $E_m,F_m$.
Since $|X\cap B_1'|\ge W'n/2$, there are at least $W'n/4-2\epsilon n$
vertices in $X\cap B_1'$ that have no neighbour in $Y'$, and in particular have no neighbour in $B_j'\cap Y'$. Since
$W'n/4-2\epsilon n>\lambda Wn$,
$B_j'\cap Y'$
$\lambda$-misses $B_1$. By the same argument it $\lambda$-misses $B_K$, and so does not $\lambda$-cover
any of $B_2\ll B_{h-1}$, since $\mathcal{B}$ is $\lambda$-concave.
In other words, for $i\in \{2 \ll h-1\}$ and $j\in \{h\ll  h+s+t-1\}$,  there are at most $\lambda Wn$ vertices in $B_i'$ with a neighbour in
$B_j'\cap Y'$; and consequently there are at most $(s+t)\lambda Wn$ vertices in $B_i'$ with a neighbour in $Y'$.
Since $|B_i'|= W'n$ and $W'n-(s+t)\lambda Wn\ge W'n/8$, this proves~(8).

\bigskip

Now there are four cases, depending on the four possible types of $U$.
First, suppose $U$ has type $(1,E)$. Let $Q'$ be the rooted tree obtained from the disjoint union of $Q(\gamma_1)$ and $S(\gamma_3)$
by adding an edge between the roots, and making the root of $Q(\gamma_1)$ the root of the new tree. Thus $Q'$
contains $Q(\gamma_1+1)$.
Each $v\in U$ is adjacent to the root, and to no other vertices, of a $\mathcal{B}$-rainbow copy
of $S(\gamma_3)$ that is contained in $G[Y']$. But from the definition
of ``anchored'', $v$ is the root of a $\mathcal{B}$-left-rainbow copy of $Q(\gamma_1)$ that is contained in
$G[Y\cup \{v\}]$. Since $Y$ is anticomplete to $Y'$, the union of these two rooted trees contains a copy of $Q(\gamma_1+1)$,
so $v$ is the root
 of a $\mathcal{B}$-left-rainbow copy of $Q(\gamma_1+1)$ that is contained in
$G[Y\cup Y'\cup \{v\}]$. Choose $B_1''\subseteq U$ of cardinality $\lceil W'n/8\rceil$;
%W(g)n/8\ge W(g+1)n,
for $h+1\le i\le h'+t-1$ choose $B_i''\subseteq B_i'$ of cardinality $\lceil W'n/8\rceil$, anticomplete to $Y'$ (this is possible by (8));
%$W(g)n-s\lambda Wn\ge W(g+1)n$
and choose $B_K''\subseteq B_K'\cap X$ of cardinality $\lceil W'n/8\rceil$ anticomplete to $Y'$ (this is possible since at least
$W'n/4 -\epsilon n$ vertices in $B_K'\cap X$ have no neighbour in
$Y'$).
%$W(g)n/4 -\epsilon n\ge W(g+1)n
By (2), $(B_i'':i\in \{1\ll h-1\}\cup \{K\})$ is a $(\gamma_1+1, \gamma_2)$-anchored minor
of $\mathcal{B}$ of width at least $W'/8$, contrary to the maximality of $\gamma_0$.

Next, suppose $U$ has type $(1,F)$. Let $Q'$ be the rooted tree obtained from the disjoint union of $Q(\gamma_1)$ and $T(\delta,\beta)$
by adding an edge between the roots, and making the root of $Q(\gamma_1)$ the root of the new tree. Again,
$Q'$ contains $Q(\gamma_1+1)$, since $\beta\ge \alpha$.
Each $v\in U$ is adjacent to the root, and to no other vertices, of a $\mathcal{B}$-rainbow copy
of $T(\delta,\beta)$ that is contained in $G[Y']$. But
$v$ is the root of a $\mathcal{B}$-left-rainbow copy of $Q(\gamma_1)$ that is contained in
$G[Y\cup \{v\}]$. The union of these two rooted trees is a copy of $Q'$,
so $v$ is the root
 of a $\mathcal{B}$-left-rainbow copy of $Q(\gamma_1+1)$ that is contained in
$G[Y\cup Y'\cup \{v\}]$. Then we obtain a contradiction as in the first case.

Next suppose $U$ has type $(K,F)$. Then similarly we obtain a $(\gamma_1, \gamma_2+1)$-anchored minor
of $\mathcal{B}$ of width at least $W'/8$, again a contradiction.

Finally, suppose $U$ has type $(K,E)$. We recall that $0\le \gamma_2\le 2\delta$. If $\gamma_2<\delta$, then as in the previous case
we obtain a $(\gamma_1, \gamma_2+1)$-anchored minor
of $\mathcal{B}$ of width at least $W'/8$, contrary to the maximality of $\gamma_0$. So we may assume that $\gamma_2\ge \delta$.
Choose $v\in U$, and let $i$ be its happiness; then $E_i$ is a copy of $S(\gamma_3)$.
Let $u$ be the root of $E_i$.  Since $v\in B_K'$, $v$ is the root of a $\mathcal{B}$-rainbow copy of $R(\gamma_2)$,
contained in $G[Y\cup \{v\}]$. But $R(\gamma_2)$ contains
$T(\delta,\alpha)$ (it even contains $T(\delta,\alpha+1)$, but we do not need that); and consequently
$v$ is the root of a $\mathcal{B}$-rainbow copy of $T(\delta,\alpha)$, contained
in $G[Y\cup \{v\}]$. The union of this tree with $E_i$, rooted at $u$, gives a $\mathcal{B}$-left-rainbow copy of $S(\gamma_3+1)$,
contrary to the choice of $\gamma_3$. This proves \ref{usemonotone}.~\bbox

\section{Producing concavity}\label{sec:proofend}

That concludes the difficult part of the paper: now we just have to write out carefully the argument sketched in section 2.

\begin{thm}\label{fixedpoint}
Let $\tau\ge 1$ be an integer, and let $0<\kappa\le 1$. 
%$c>0$ satisfy $$\log(1/c)=m^{\nu+1}\log (1/\mu_0).$$ (Logarithms are to base two.)
Let $\mathcal{B}=(B_1\ll B_K)$ be a blockade in a graph.
Then there is an equicardinal $(\kappa,\tau)$-support-invariant contraction $\mathcal{B}'=(B_1'\ll B_K')$ of $\mathcal{B}$
with width at least $\kappa^{2^K\tau^\tau}$ times the width of $\mathcal{B}$.
\end{thm}
\Proof
If $\tau\ge 1$ is an integer, we define the {\em $\tau$-cost}
of a blockade $\mathcal{B}'$ to be the sum of the cardinalities of the traces of all nonisomorphic ordered trees $J$
with at most $\tau$ vertices. Since there are only at most $\tau^\tau$ nonisomorphic ordered trees with at most $\tau$ vertices,
and the trace of each relative to $\mathcal{B}'$ has cardinality at most $2^K$, 
the $\tau$-cost of $\mathcal{B}'$ is at most $2^K\tau^\tau$. Let the width of $\mathcal{B}$ be $W$. 
Choose an integer $t\le 2^K\tau^\tau$, maximum such that there is an equicardinal contraction $\mathcal{B}'=(B_1'\ll B_K')$ of $\mathcal{B}$
with width at least $\kappa^tW$ and $\tau$-cost at most  $2^K\tau^\tau -t$. We claim that $\mathcal{B}'$ is 
$(\kappa,\tau)$-support-invariant. Suppose not; then there is an ordered tree $J$ with $|J|\le \tau$,
and a 
contraction $\mathcal{B}''=(B_i'':i\in I)$ of $\mathcal{B}'$ of width at least $\kappa$ times the width of $\mathcal{B}'$,
and the trace of $J$ relative to $\mathcal{B}''$ is different from the trace of $J$ relative to $\mathcal{B}'$. But the first
trace is a subset of the second, so it is a proper subset.
For every other ordered tree $J'$, its trace relative to $\mathcal{B}''$
is a subset of its trace relative to $\mathcal{B}'$; and so the $\tau$-cost of $\mathcal{B}''$ is strictly less than the $\tau$-cost of 
$\mathcal{B}'$. But this contradicts the maximality of $t$, and so proves \ref{fixedpoint}.~\bbox

By iterated applications of Ramsey's theorem for uniform hypergraphs (applying it
once for each ordered tree $J$ with $|J|\leq \tau$ and dropping to a suitable
sub-blockade each time)
we deduce:

\begin{thm}\label{getuniform}
Let $k\ge 0$ and $\tau\ge 1$ be integers; then there exists an integer $K\ge 0$ with the following property. Let 
$\mathcal{B}=(B_1\ll B_K)$ be a blockade
in a graph. Then $\mathcal{B}$ 
has a sub-blockade of length $k$ which is $\tau$-support-uniform.
\end{thm}

Combining \ref{fixedpoint} and \ref{getuniform}, we obtain:

\begin{thm}\label{getminor}
Let $k\ge 0$ and $\tau\ge 1$ be integers, and $0<\kappa\le 1$; then there exist an integer $K$ with the following property. 
Let $\mathcal{B}=(B_1\ll B_K)$ be a blockade in a graph, and let $W$ be its width.
Then there is an equicardinal minor $\mathcal{B}'$ of $\mathcal{B}$, with 
length $k$ and width at least $\kappa^{2^K\tau^\tau}W$, such that 
$\mathcal{B}'$ is $\tau$-support-uniform and $(\kappa,\tau)$-support-invariant.
\end{thm}
\Proof
Let $K$ satisfy \ref{getuniform};
then we claim it satisfies \ref{getminor}. Let 
$\mathcal{B}=(B_1\ll B_K)$ be a blockade of width $W$ in a graph. By \ref{fixedpoint} 
there is an equicardinal $(\kappa,\tau)$-support-invariant contraction $\mathcal{B}'$ of $\mathcal{B}$, with width at least 
$\kappa^{2^K\tau^\tau}W$.
By \ref{getuniform} applied to $\mathcal{B}'$, the result follows, since being $(\kappa,\tau)$-support-invariant is inherited by
sub-blockades. This proves \ref{getminor}.~\bbox

Now we prove \ref{rainbow}, which we restate:
\begin{thm}\label{rainbow2}
For every tree $T$, there exist $d>0$ and an integer $K$, such that,
for every graph $G$ with a blockade $\mathcal{B}$ of length at least $K$, if $G$ is $\frac{W}{d|G|}$-coherent where $W$ is the width
of $\mathcal{B}$, then
there is a $\mathcal{B}$-rainbow copy of $T$ in $G$.
\end{thm}
\Proof
We proceed by induction on $|T|$, and may assume that $|T|\ge 2$. Choose $\delta\ge 2$ and $\eta\ge 0$ such that $T$ is a 
subtree of $T(\delta,\eta)$. Let $\lambda = 2^{-9\delta}\delta^{-1-\eta}$, and $\kappa=\lambda/2$.
Let $r= \lceil (|T|-1)/\kappa\rceil$.
Let $v$ be a vertex of $T$ with degree one, and let $u$
be its neighbour. Let $T'=T\setminus \{v\}$; from the inductive hypothesis, there exist $K', d'$ satisfying the theorem with $T$
replaced by $T'$. By increasing $K'$, we may assume that $K'\ge 6r\delta^{\eta+2}$, and $K'$ is a multiple of $r$. 
Let $k=K'/r$.
Let $\tau=\delta^{\eta+1}$.
Let $K$ satisfy \ref{getminor} with $k$ replaced by $K'$.
Let 
$$d=\kappa^{-2^K\tau^\tau}\max(d', 2^{9\delta}/r).$$
%$\kappa^{\tau^\tau2^K}/d'\ge 1/d$
%$r\kappa^{\tau^\tau2^K}\ge 2^{9\delta}/d$
We claim that $K,d$ satisfy the theorem. Let $\mathcal{B}$
be a blockade in a graph $G$, of length $K$ and width $W$, such that $G$ is $\frac{W}{d|G|}$-coherent; and let
$\vare=W/(d|G|)$.
We assume (for a contradiction) that there
is no $\mathcal{B}$-rainbow copy of $T$. By \ref{getminor}, there is an equicardinal minor $\mathcal{B}'$
of $\mathcal{B}$ of length $K'$ and width at least $\kappa^{2^K\tau^\tau}W$, such that
$\mathcal{B}'$ is $\tau$-support-uniform and $(\kappa,\tau)$-support-invariant. 
Let $\mathcal{B}'=(B_1'\ll B_{K'}')$, and let its width be $w$.

From the inductive hypothesis, there is a $\mathcal{B}'$-rainbow copy of $T'$, since
$$\kappa^{2^K\tau^\tau}W\ge \kappa^{2^K\tau^\tau}(\vare d |G|)\ge d'\vare|G|.$$
%k\ge K'$ 
%$\kappa^{\tau^\tau2^K}W/d'\ge W/d$
Hence there is an ordered graph $J$, obtained from $T'$ by ordering its vertices, with nonempty trace relative to $\mathcal{B}'$.
Let $u$ be the $j$th vertex in the ordering of the vertices of $J$, and let $t=|J|=|T|-1$. Thus $t\le \tau$.
\\
\\
(1) {\em There do not exist $1\le r_1<\cdots< r_t\le K'$, such that for some $X\subseteq V(\mathcal{B}')$, $X$ is disjoint from
$B_{r_h}'$ for $1\le h\le t$, and 
$X$ $\kappa$-covers $B_{r_j}'$, and $X$ $\kappa$-misses $B_{r_h}'$ for all $h\in \{1\ll t\}\setminus \{j\}$.}
\\
\\
Suppose that such $X$ and $r_1\ll r_t$ exist. Let $B_{r_j}''$ be the set of vertices in $B_{r_j'}$ that have a neighbour in $X$;
for $1\le h\le t$ with $h\ne j$, let $B_{r_h}''$ be the set of vertices in $B_{r_h'}$ that do not have a neighbour in $X$; and for 
$h\in \{1\ll K'\}$ with $h\ne r_1\ll r_t$, let $B_h''=B_h'$. Then $(B_1''\ll B_{K'}'')$ has width at least
$\kappa w$. 
Since  $\mathcal{B}'$ is $\tau$-support-uniform, and $(\kappa,\tau)$-support-invariant,
and the trace of $J$  relative to $\mathcal{B}'$
is nonempty, it follows that there is a $\mathcal{B}'$-rainbow induced subgraph $H$ of $G$, with $\mathcal{B}'$-ordering
isomorphic to $J$, where for $1\le h\le t$ the $h$th
vertex of $H$ belongs to $B_{r_h}''$. In particular, there is an isomorphism from $T'$ to $H$ mapping $u$
to the vertex, $u'$ say, of $H$ in $B_{r_j}''$. Choose $v'\in X$ adjacent to $u'$; such a vertex exists since 
$X$ covers $B_{r_j}''$. But then $v'$ has no other neighbour in $V(H)$, since $X$ misses $B_{r_h}''$ 
for all $1\le h\le t$ with $h\ne j$; and $v'\in B_i'$ for some $i\in \{1\ll K'\}$ with $i\ne r_1\ll r_t$. 
Thus adding $v'$ to $H$ gives a $\mathcal{B}$-rainbow copy of $T$, a contradiction. This proves (1).

\bigskip
We recall that $K'=rk$.
For $1\le h\le k$, let
%$k'=rk
$C_h$ be the union of the sets $B_i'$ for all $i$ with $r(h-1)<i\le rh$.
Then $\mathcal{C}=(C_1\ll C_k)$ is a blockade, of width $rw\ge r\kappa^{2^K\tau^\tau}W$. 
\\
\\
(2) {\em $\mathcal{C}$ is $2\kappa$-concave.}
\\
\\
Let $1\le h_1<h_2<h_3\le k$, and suppose that there is a set $X\subseteq V(\mathcal{C})$ such that
$X$ $\lambda$-covers $C_{h_2}$ and $\lambda$-misses $C_{h_1}$ and $C_{h_3}$, and $X$ is disjoint from $C_{h_1}\cup C_{h_2}\cup C_{h_3}$.
Since $|X|$ $\lambda$-covers $C_{h_2}$, there are at least $\lambda rw$ vertices in $C_{h_2}$ with a neighbour in $X$, and so there
exists $r_j$ with $r(h_2-1)< r_j\le rh_2$ such that at least $\lambda w$ vertices in $B_{r_j}'$ have a neighbour in $X$; and so
$X$ $\lambda$-covers $B_{r_j}'$, and so $X$ $\kappa$-covers $B_{r_j}'$ . Hence by (1), either $X$ $\kappa$-misses $B_i'$ for fewer than $j-1$ (and hence fewer than $t$) 
%\kappa\le \lambda$
values of 
$i$ with $r(h_1-1)<i\le rh_1$, or $X$ $\kappa$-misses $B_i'$ for fewer than $t-j$ (and hence fewer than $t$) values of 
$i$ with $r(h_3-1)<i\ll rh_3$; and from the symmetry we may assume the former. Consequently 
there are fewer than 
$tw+(r-t)\kappa w$ vertices in $C_{h_1-1}$ with  no neighbour in $X$. Since $X$ $\lambda$-misses $C_{h_1}$, it follows that
$$tw+(r-t)\kappa w> \lambda|C_{h_1}|=2\kappa rw,$$
and so 
$t(1-\kappa)> \kappa r,$
contrary to the choice of $r$. This proves (2).
%r\ge t/\kappa
%$\lambda=2\kappa
\\
\\
(3) {\em $\mathcal{C}$ is $\tau$-support-uniform and $(\kappa,\tau)$-support-invariant.}
\\
\\
The first statement is clear, and we only need prove the second. Let $S$ be an ordered tree, with $|S|\le \tau$,
such that there is a $\mathcal{C}$-rainbow copy of $S$. Let $|S|=s$ say.
Let $\mathcal{C}'=(C_{h_1}'\ll C_{h_s}')$ be a contraction of a sub-blockade of $\mathcal{C}$, of width at least $\kappa rw$, and
of length $s$. We must show that there is a $\mathcal{C}'$-rainbow copy of $S$. To simplify notation
we assume without loss of generality that $h_i=i$ for $1\le i\le s$. Now for $1\le i\le s$, since $|C_i'|\ge \kappa rw$,
there exists $g_i$ with $r(i-1)< g_i\le ri$ such that $|C_i'\cap B'_{g_i}|\ge \kappa w$.
Since $\mathcal{B}'$ is $\tau$-support-uniform and $(\kappa,\tau)$-support-invariant, and $|J|\le \tau$, it follows
that there is a copy of $S$ that is $(C_1'\cap B'_{g_1}\ll C_s'\cap B'_{g_s})$-rainbow, and hence $(C_1'\ll C_s')$-rainbow.
This proves (3).

\bigskip

From (2) and (3) we see that $\mathcal{C}$ is $\lambda$-concave, $\tau$-support-uniform and $(\kappa,\tau)$-support-invariant, and
therefore $(2^{-9\delta},\tau)$-support-invariant, since $\kappa\le 2^{-9\delta}$.
Its length is $k\ge 6\delta^{\eta+2}$, and its width is at least $r\kappa^{2^K\tau^\tau}W\ge 2^{9\delta}\epsilon|G|$.
%$r\kappa^{\tau^\tau2^K}W\ge 2^{9\delta}\epsilon$
From \ref{usemonotone}, there is a $\mathcal{C}$-rainbow copy of $T(\delta,\eta)$ and hence of $T$. Consequently this copy of $T$ is also
$\mathcal{B}'$-rainbow, and therefore $\mathcal{B}$-rainbow. This proves \ref{rainbow2}.~\bbox

Finally we can prove \ref{mainthm}, which we restate:

\begin{thm}\label{mainthm2}
For every tree $T$ there exists $\epsilon>0$ such that every $\epsilon$-coherent graph contains $T$.
\end{thm}
\Proof 
Let $T$ be a tree, and choose $\delta\ge 2$ and $\eta\ge 0$ such that $T(\delta,\eta)$ contains $T$.
Let $K,d$ satisfy \ref{rainbow2}; we may assume by increasing $d$ that $d\ge 1$. Choose $\epsilon>0$ such that $2K\epsilon d\le 1$.
We claim that every $\epsilon$-coherent graph contains $T$.
Let $G$ be $\epsilon$-coherent; it follows that $|G|\ge \epsilon^{-1}\ge 2Kd\ge 2K$. Hence
$|G|/K\ge \lceil |G|/(2K)\rceil$, and so we may choose $K$ subsets of $V(G)$, pairwise disjoint and
each of cardinality at least $|G|/(2K) \ge \epsilon d |G|$. These sets, in any order, form a blockade of length $K$ 
and width at least $\epsilon d|G|$,
and so by \ref{rainbow2}, $G$ contains $T$. 
This proves \ref{mainthm2}.~\bbox

\section{Remarks}

There are some final points we would like to make. First, while the results of \cite{sparse,cats} concerned
$\epsilon$-coherent graphs, they were capable of generalization in the natural way
to $\epsilon$-coherent ``massed graphs'', graphs in which each subset $X\subseteq V(G)$ had a mass $\mu(X)$,
where $\mu$ was increasing and subadditive (and also satisfied a nontriviality condition); such as, for instance, the
function $\mu(X)=\chi(G[X])/\chi(G)$, where $\chi$ denotes chromatic number. The proof of \ref{mainthm} does not seem to extend
to massed graphs; for instance, the method in the proof of \ref{usemonotone} of pulling out ``parallel'' rainbow copies of 
a graph, relies on the fact that we are removing the same number of vertices from each block. 

Second, the following was proposed as a conjecture in \cite{sparse}, and we have now found a proof (closely related to the proof in this paper),
which will appear in~\cite{pure4}:

\begin{thm}\label{bipconj}
For every forest $T$ there exists $\epsilon>0$ with the following property.
Let $G$ be a $T$-free bipartite graph with bipartition $(A,B)$, where $|A|=|B|=n$. Then either
some vertex has degree at least $\epsilon n$, or
there is an anticomplete pair of subsets $A'\subseteq A$ and $B'\subseteq B$ with $|A'|, |B'|\ge \epsilon n$.
\end{thm}

Third, here is a nice question: for which tournaments $H$ does there exist $\epsilon>0$ such that in every tournament $G$ with $|G|>1$
not containing $H$ as a sub-tournament, there are two disjoint subsets $A,B$ of $V(G)$ where $A$ is complete to $B$ and 
$|A|,|B|\ge \epsilon |G|$? (If so, we say $H$ has the ``strong EH-property''.) One can show that if $H$ 
has the strong EH-property, then
\begin{itemize}
\item $V(H)$ can be ordered as $\{v_1\ll v_n\}$ such that the backedge digraph (the digraph formed by the pairs $v_iv_j$ where $v_i$
is adjacent from $v_j$ in $G$) is transitive;
\item $V(H)$ can be ordered such that the backedge digraph has no induced outdirected 3-star;
\item $V(H)$ can be ordered such that the backedge digraph has no induced indirected 3-star; and
\item $V(H)$ can be ordered such that the backedge graph (the graph underlying the backedge digraph) is a forest.
\end{itemize}
and some other similar conditions.
But such tournaments exist; for instance, the eulerian orientation of $K_5$ is such a tournament, and recently Berger, Choromanski, Chudnovsky and Zerbib~\cite{berger} 
proved that it has the strong EH-property. The seven-vertex Paley tournament does not  have the strong EH-property, but if we delete 
one vertex from it, we obtain a six-vertex tournament that might have the property, and does satisfy the bullets above; this is one of the two smallest tournaments that are currently undecided.

Fourth, for a graph $H$, define $d(H)$ to be the minimum of $(|J|-1)/|E(J)|$ over all induced subgraphs $J$ of $H$
that have at least one edge. Thus $d(H)\le 1$, with equality if and only if $H$ is a forest. It is tempting to conjecture 
that for all $H$, there exists $\epsilon>0$
such that in every $H$-free graph $G$ with $|G|>1$ and maximum degree less than $\epsilon |G|$, 
there are two disjoint subsets $A,B\subseteq V(G)$,
anticomplete, with $|A|,|B|\ge \epsilon |G|^{d(H)}$. When $d(H)=1$ this is our theorem, and one can
modify Erd\H{o}s' random graph construction to show that the bound would be sharp for all $d(H)$. Unfortunately it is false;
for instance, when $H=K_3$, there is a sparse $H$-free graph $G$ with $n$ vertices for large $n$,
in which every anticomplete pair of sets $A,B$ satisfies $\min(|A|,|B|)\le O(n^{1/2}\polylog(n))$. 
Nevertheless, $d(H)$ gives a guide: we will show in a later paper~\cite{pure8} that for all $c$ with $0<c<1$, if $H$ is a graph with 
$1-d(H)$ sufficiently
small, then there exists 
$\epsilon>0$ such that in every $H$-free graph $G$ with $|G|>1$ and maximum degree less than $\epsilon |G|$, there are two disjoint anticomplete subsets $A,B\subseteq V(G)$
with $|A|,|B|\ge \epsilon |G|^{c}$.

Fifth, what about ordered graphs? Let $H$ be an ordered forest, and let us look at sparse ordered 
graphs
$G$ that do not contain $H$ as an ordered induced subgraph, in the natural sense. Is it still true that there are two
disjoint anticomplete subsets of linear cardinality? No, it is not: but we will show in a later paper~\cite{pure6} that there are 
two anticomplete sets both of cardinality at least $|G|^{1-o(1)}$. Also, in~\cite{pure7} we will show a similar theorem modifying \ref{bipconj},
for ordered
bipartite graphs (placing an ordering on both sets of the bipartition, not on their union).

Finally, as with many Ramsey-type theorems, there is a multicolouring version of our result. Take a complete graph, and partition
its edge-set into $k$ sets; and let $G_i$ be the subgraph with edge-set the $i$th of these sets (and all the vertices). 
We call $(G_1\ll G_k)$ a {\em $k$-multicolouring}. Then the following holds, generalizing \ref{forestsymm}:
\begin{thm}\label{multicolour}
For all $k\ge 1$ and every forest $H$ there exists $\epsilon>0$, such that if $(G_1\ll G_k)$ is a $k$-multicolouring of a complete 
graph $K_n$ with at least two vertices, then for some $i\in \{1\ll k\}$, either $G_i$ contains $H$ as an induced subgraph, 
or there are two disjoint subsets 
$X,Y\subseteq V(G_i)$, with $|X|,|Y|\ge \epsilon n$, anticomplete in $G_i$.
\end{thm}
\Proof (Sketch.) The proof is an easy corollary of \ref{mainthm}. Choose $\epsilon'>0$ such that \ref{mainthm} holds, and choose $c>0$ such that $kc \le \epsilon'$.
A straightforward modification of the proof of the theorem of \cite{rodl} shows that
there exists $\delta>0$ (independent of $n$ and $G_1\ll G_k$), such that if no $G_i$ contains $H$
as an induced subgraph, then there is a subset $A$ of the vertex set of $K_n$, with $|A|\ge \delta n$, such that 
$|E(G_i[A])|\le c|A|(|A|-1)/2$
for all values of $i\in \{1\ll k\}$ except one, say all except $i=k$. Let $\epsilon=\epsilon'\delta$. 
By removing vertices of degree at least $kc\delta n$ in one of 
$G_1\ll G_{k-1}$, we deduce that there exists $B\subseteq A$ with $|B|\ge (\delta/k) n$, such that 
for $1\le i\le k-1$, every vertex of $G_i[B]$ has degree at most $kc\delta n\le \epsilon'|B|$ in $G_i[B]$. By \ref{mainthm} applied
to $G_1[B]$,
there are two disjoint subsets $X,Y\subseteq B$ with $|X|,|Y|\ge \epsilon'|B|\ge \epsilon n$, anticomplete 
in $G_1[B]$, as required.~\bbox

There is in fact a stronger result:
\begin{thm}\label{multicolour2}
For all $k\ge 2$ and every forest $H$ there exists $\epsilon>0$, such that if $(G_1\ll G_k)$ is a $k$-multicolouring of a complete
graph $K_n$ with at least two vertices, then for some distinct $i,j\in \{1\ll k\}$, either 
\begin{itemize}
\item there is a subset $X\subseteq V(K_n)$ such that every edge of $G[X]$ belongs to $G_i\cup G_j$, and $G_i[X]$
is isomorphic to $H$; or
\item there are two disjoint subsets
$X,Y\subseteq V(G_i)$, with $|X|,|Y|\ge \epsilon n$, complete in $G_j$.
\end{itemize}
\end{thm}
We do not know how to deduce this as a consequence of \ref{mainthm}, and it seems necessary to modify the proof of
\ref{mainthm} in several places; all straightforward, but too many to sketch here, and we omit further details.


\begin{thebibliography}{99}
\bibitem{APPRS}
N. Alon, J. Pach, R. Pinchasi, R. Radoi\v ci\'c and M. Sharir, ''Crossing patterns of semi-algebraic
sets'', {\em J. Combinatorial Theory, Ser. A}, {\bf 111} (2005), 31--326.
\bibitem{berger} E. Berger, K. Choromanski, M. Chudnovsky and S. Zerbib, ``Tournaments and the strong Erd\H{o}s-Hajnal property'', 
{\tt arXiv:2002.07248}.
\bibitem{bonamy} M. Bonamy, N. Bousquet and S. Thomass\'e, ``The Erd\H{o}s-Hajnal conjecture for long holes and antiholes'',
{\em SIAM J. Discrete Math.} {\bf 30} (2015), 1159--1164.
\bibitem{lagoutte} N. Bousquet, A. Lagoutte, and S. Thomass\'e, ``The Erd\H{o}s-Hajnal conjecture for paths and
antipaths'', {\em J. Combinatorial Theory, Ser. B}, {\bf 113} (2015), 261--264.
\bibitem{hooks} K. Choromanski, D. Falik, A. Liebenau, V. Patel, and M. Pilipczuk, ``Excluding hooks and their complements'',
{\em Electronic J. Combinatorics} {\bf 25} \#P3.27, {\tt arXiv:1508.00634}.
\bibitem{sparse} M. Chudnovsky, A. Scott, P. Seymour and S. Spirkl, ``Pure pairs. II. Excluding all subdivisions of a graph'', submitted for publication, {\tt arXiv:1804.01060}.
\bibitem{girth} P. Erd\H{o}s, ``Graph theory and probability'', {\em Canadian J. Math.} {\bf 11} (1959), 34--38.
\bibitem{EH0} P. Erd\H{o}s and A. Hajnal, ``On spanned subgraphs of graphs'',
{\em Graphentheorie und Ihre Anwendungen} (Oberhof, 1977), \verb++{www.renyi.hu/\raisebox{-1ex}{\textasciitilde}p\_erdos/1977-19.pdf}.
\bibitem{EH}  P. Erd\H{o}s and A. Hajnal, ``Ramsey-type theorems'',
{\em  Discrete Applied Math.} {\bf 25} (1989), 37--52.
\bibitem{fp}
J. Fox and J. Pach, ``Erd\H{o}s-Hajnal-type results on intersection patterns of geometric objects'', in
{\em Horizon of Combinatorics} (G.O.H. Katona et al., eds.), Bolyai Society Studies in Mathematics,
Springer, 79--103, 2008.
\bibitem{cats} A. Liebenau, M. Pilipczuk, P. Seymour and S. Spirkl, ``Caterpillars in Erd\H{o}s-Hajnal'',
{\em J. Combinatorial Theory, Ser. B}, {\bf 136} (2019), 33--43, {\tt arXiv:1810.00811}.
\bibitem{rodl} V. R\"odl, ``On universality of graphs with uniformly distributed edges'',
{\em Discrete Math.} {\bf 59} (1986), 125--134.
\bibitem{pure4} A. Scott, P. Seymour and S. Spirkl, ``Pure pairs. IV. Trees in bipartite graphs'', in preparation.
\bibitem{pure6} A. Scott, P. Seymour and S. Spirkl, ``Pure pairs. VI. Excluding an ordered tree'', in preparation.
\bibitem{pure7} A. Scott, P. Seymour and S. Spirkl, ``Pure pairs. VII. Homogeneous submatrices in a 0/1-matrix with a forbidden submatrix'', in preparation.
\bibitem{pure8} A. Scott, P. Seymour and S. Spirkl, ``Pure pairs. VIII. Excluding a sparse graph'', in preparation.



\end{thebibliography}
\end{document}